\documentstyle[12pt]{amsart}
\include{mssymb}
\newtheorem{prop}{Proposition}
\newtheorem{def-prop}{Definition-Proposition}
\newtheorem{thm}{Theorem}
\newtheorem{cor}{Corollary}

\newtheorem{lemma}{Lemma}
\newtheorem{defn}{Definition}
\def\endproof{$\Box$ \medskip}
\def\1{1}
\def\P{{\Bbb P}}
\def\1{{\Bbb I}}

\def\Z{{\Bbb Z}}

\def\G{{\Bbb G}}

\def\Q{{\Bbb Q}}

\def\B{{\cal B}}
\def\Ad{\mbox{Ad }}
\def\g{{\frak g}}
\def\b{{\frak b}}

\begin{document}
\title[Localization and the Bott residue formula]
{Localization in equivariant intersection theory and the Bott residue
formula}
\author{Dan Edidin}
\address{Department of Mathematics\\ University of Missouri\\
Columbia MO 65211}
\author{William Graham}
\address{School of Mathematics\\ Institute for Advanced Study\\
Princeton NJ 08540}
\thanks{The first author was partially by an NSF post-doctoral
fellowship and the M.U. research board. The second author
was partially supported by NSF grant DMS 9304580.}
\date{}
\maketitle
\section{Introduction}

The purpose of this paper is to prove the localization theorem for
torus actions in equivariant intersection theory.  Using the theorem
we give another proof of the Bott residue formula for Chern numbers of
bundles on smooth complete varieties. In addition, our
techniques allow us to obtain residue formulas for bundles on a
certain class of singular schemes which admit torus actions. This
class is rather special, but it includes some interesting examples
such as complete intersections (cf. \cite{BFQ}) and Schubert varieties.

Let $T$ be a split torus acting on a scheme $X$.  The $T$-equivariant
Chow groups of $X$ are a module over $R_T = Sym(\hat{T})$, where
$\hat{T}$ is the character group of $T$.  The localization theorem
states that up to $R_T$-torsion, the equivariant Chow groups of the
fixed locus $X^T$ are isomorphic to those of $X$.  Such a theorem is a
hallmark of any equivariant theory. The earliest version (for
equivariant cohomology) is due to Borel \cite{Borel}.  Subsequently
$K$-theory versions were proved by Segal \cite{Segal} (in topological
$K$-theory), Quart \cite{Quart} (for actions of a cyclic group), and
Thomason  \cite{Thomason} (for algebraic $K$-theory \cite{Thomason}).

For equivariant Chow groups, the localization isomorphism is given by
the equivariant pushforward $i_*$ induced by the inclusion of $X^T$ to
$X$.  An interesting aspect of this theory is that the push-forward is
naturally defined on the level of cycles, even in the singular case.
The closest topological analogue of this is equivariant Borel-Moore
homology (see \cite{E-G} for a definition), and a similar proof
establishes localization in that theory.

For smooth spaces, the inverse to the equivariant push-forward can be
written explicitly.  It was realized independently by several authors
(\cite{I-N}, \cite{A-B}, \cite{B-V}) that for compact spaces, the
formula for the inverse implies the Bott residue formula.  In this
paper, we prove the Bott residue formula for actions of split tori on
smooth complete varieties defined over an arbitrary field, also by
computing $(i_*)^{-1}$ explicitly.  Bott's residue formula has been
applied recently in enumerative geometry (cf. \cite{E-S}, \cite{K})
and there was interest in a purely intersection-theoretic proof.
Another application of the explicit formula for $(i_*)^{-1}$ is given
in \cite{E-G2}, where we prove (following Lerman \cite{L}) a residue
formula due to Kalkman.

An obvious problem, which should have applications to enumerative
geometry (see e.g. \cite{K}), is to extend the Bott residue formula to
complete singular schemes.  Such a formula can be derived when we have
an explicit description of $(i_*)^{-1}[X]_T$, where $[X]_T$ denotes
the equivariant fundamental class of the whole scheme, as follows: Let
$n = \mbox{dim }X$.  If $[X]_T = i_*\alpha$ and $p(E)$ is a polynomial
of weighted degree $n$ in Chern classes of equivariant vector bundles
on $X$, then $\mbox{deg }(p(E) \cap [X])$ can be calculated as the
residue of $\pi_*(i^*(p(E)) \cap \beta)$ where $\pi$ is the
equviariant projection from the $X$ (or the fixed locus) to a point.
This approach does not work for equivariant cohomology, because when
$X$ is singular there is no pushforward from $H^*_G(X) \rightarrow
H^*_G(M)$.  However, in $K$-theory, where such pushforwards exist,
similar ideas were used by \cite{BFQ} to obtain Lefschetz-Riemann-Roch
formulas for the action of an automorphism of finite order.

The problem of computing $(i_*)^{-1}$ is difficult, but we can do it
in a certain class of singular examples, in particular, if there is an
equivariant embedding $X \stackrel{f}\hookrightarrow M$
into a smooth variety, and every
component of $X^T$ is a component of $M^T$. This condition is
satisfied if $X \subset \P^r$ is an invariant subvariety where $T$
acts linearly with distinct weights (and thus isolated fixed points)
or if $X$ is a Schubert variety in $G/B$.  In this context we give a
formula (Proposition \ref{xxxsing}) formula for $(i_*)^{-1}\alpha$ in
terms of $f_*\alpha \in A_*^T(M)$. The case of Schubert
of varieties is worked out in detail in Section \ref{schubs}.

As a consequence it is possible to compute Chern numbers of bundles on
$X$ provided we know $f_*[X]_T \in A_*^T(M)$. Thus for example, if $X$
is a $T$-invariant projective variety and $T$ acts linearly with
distinct weights on $\P^n$, then we can calculate Chern numbers,
provided we know the equivariant fundamental class of $X$.  Rather
than write down a general formula, we illustrate this with an
example: in Section \ref{singex} we use a residue calculation to show
that
$$
\int_Q c_1(\pi^* T_{\P^2}) c_1(f^*T_{\P^3}) = 24
$$
where $Q \stackrel{f} \hookrightarrow \P^3$ is a (singular) quadric cone,
and
$\pi: Q \rightarrow \P^2$ is the projection from a point not on $Q$.
The methods can be applied in other examples.

{\bf Acknowledgements.} We thank
Steven Kleiman for suggesting the problem of giving an
algebraic proof of the Bott residue formula. We are also
grateful to Robert Laterveer for discussions of Gillet's
higher Chow groups.

\section{Review of equivariant Chow groups}

In this section we review some of the equivariant intersection
theory developed in \cite{E-G}. The key to
the theory is the definition of equivariant
Chow groups for actions of linear algebraic groups.
All schemes are assumed to be of finite type
defined over a field of arbitrary characteristic.

Let $G$ be a $g$-dimensional group, $X$ an $n$-dimensional scheme and
$V$ a representation of $G$ of dimension $l$. Assume that there is an
open set $U \subset V$ such that a principal bundle quotient $U
\rightarrow U/G$ exist, and that $V-U$ has codimension more than
$n-i$. Thus the group $G$ acts freely on the product $X \times U$. The
group $G$ acts freely on $X \times U$, and if any one of a number of
mild hypotheses is satisfied then a quotient scheme $X_G = (X \times
U)/G$ exists (\cite{E-G}).  In particular, if $G$ is special -- for
example, if $G$ is a split torus, the case of interest in this paper
--  a quotient scheme $X_G$ exists.

\begin{defn}
Set $A_i^G(X) = A_{i+l-g}(X_G)$, where $A_*$ is the usual Chow group.
This definition is independent of the choice of $V$ and $U$ as long as
$V-U$ has sufficiently high codimension.
\end{defn}

{\bf Remark:} Because $X \times U \rightarrow X \times^G U$ is
a principal $G$-bundle, cycles on $X \times^G U$ exactly
correspond to $G$-invariant cycles on $X \times U$. Since
we only consider cycles of codimension smaller
than the dimension of $X \times (V-U)$, we may in fact
view these as $G$-invariant cycles on $X \times V$. Thus
every class in $A_i^G(X)$ is represented by a cycle in
$Z_{i+l}(X \times V)^G$, where $Z_*(X \times V)^G$ indicates
the group of cycles generated by invariant subvarieties.
Conversely, any cycle in $Z_{i+l}(X \times V)^G$ determines
an equivariant class in $A_i^G(X)$.

\medskip

The properties of equivariant intersection Chow groups include the following.

(1) Functoriality for equivariant maps:
proper pushforward, flat pullback, l.c.i pullback, etc.

(2) Chern classes of equivariant bundles operate on  equivariant
Chow groups.

(3) If $X$ is smooth of dimension $n$, then we denote $A_{n-i}^G(X)$ as
$A^i_G(X)$. In this case there is an intersection
product $A^i_G(X) \times A^j_G(X) \rightarrow A_G^{i+j}(X)$,
so the groups $\oplus_0^\infty A^i_G(X)$ form a graded ring which
we call the equivariant Chow ring. Unlike, the ordinary case
$A^i_G(X)$ can be non-zero for any $i \geq 0$.
(The existence of an intersection product follows from (1),
since the diagonal $X \hookrightarrow X \times X$
is an equivariant regular embedding when $X$ is smooth.)

(4) Of particular use for this paper is the equivariant self-intersection
formula. If $Y \stackrel{i} \hookrightarrow X$ is a regularly embedded
invariant subvariety of codimension $d$,
then
$$i^*i_*(\alpha) = c_d^G(N_YX) \cap \alpha$$ for any
$\alpha \in A_*^G(Y)$.

\subsection{Equivariant higher Chow groups}
Let $Y$ be a scheme. Denote by $A_i(Y,j)$ the higher
Chow groups of Bloch \cite{Bl} (indexed by dimension)
or the groups $CH_{i,i-j}(X)$ defined in \cite[Section 8]{Gillet}.
Both theories agree with ordinary Chow groups when $j=0$,
and both extend the localization short exact sequence for
ordinary Chow groups. However, in the case of Bloch's
Chow groups the localization exact sequence has only
been proved for quasi-projective varieties. The advantage
of his groups is that are naturally defined in terms
of cycles on $X \times \Delta^j$ (where $\Delta^j$ is an algebraic
$j$-simplex) and are rationally isomorphic to higher $K$-theory.

Both these theories can be extended to the equivariant
setting. We define the higher Chow groups $A_i^G(X,j)$
as $A_{i+l-g}(X_G,j)$ for an appropriate mixed space
$X_G$. Because of the quasi-projective hypothesis
in Bloch's work, Bloch's equivariant higher Chow groups
are only defined for (quasi)-projective varieties with
linearized actions. However, Gillet's are defined
for arbitrary schemes with a $G$-action. We
will use two properties of the higher equivariant
theories.\\

(a) If $E \rightarrow X$ an equivariant vector bundle, then
the equivariant Chern classes $c_i^G(E)$ operate
on $A_*^G(X,j)$.\\

(b) If $U \subset X$ is an invariant open set, then there
is a long exact sequence
$$ \ldots \rightarrow A_i^G(U,1) \rightarrow A_i^G(X-U) \rightarrow
A_i^G(X) \rightarrow A_i^G(U) \rightarrow 0.$$

\section{Localization}
In this section we prove the main theorem
of the paper,
the localization theorem for equivariant Chow groups.
For the remainder of the paper, all tori are assumed
to be split, and the coefficients off all Chow groups are rational.

Let $R_T$ denote the $T$-equivariant Chow ring of a point, and let
$\hat{T}$ be the character group of $T$.
\begin{prop} \cite[Lemma 4]{EGT}
$R_T = Sym(\hat{T})\simeq \Q[t_1, \ldots , t_n]$.
where $n$ is the rank of $T$.
\endproof \end{prop}

{\bf Remark.} The identification $R_T = Sym(\hat{T})$
is given explicitly as follows. If $\lambda \in \hat{T}$ is
a character, let $k_\lambda$ be the corresponding 1-dimensional
representation and let $L_\lambda$ denote the line bundle
$U \times^T k_\lambda \rightarrow U/T$. The map
$\hat{T} \rightarrow R^1_T$ given by $\lambda \mapsto c_1(L_{\lambda})$
extends to a ring isomorphism $Sym(\hat{T}) \rightarrow R_T$.

\medskip

\begin{prop}
If $T$ acts trivially on $X$, then $A_*^T(X) = A_*(X) \otimes
R_T$.
\end{prop}
Proof. If the action is trivial then $(U \times X)/T= U/T \times X$.
The spaces $U/T$ can be taken to be products
of projective spaces, so $A_*(U/T \times X) = A_*(X) \otimes A_*(U/T)$.
\endproof

\medskip

If $T \stackrel{f} \rightarrow  S$ is a homomorphism of tori, there
is a pullback $\hat{S} \stackrel{f^*} \rightarrow \hat{T}$. This
extends to a ring homomorphism $Sym(\hat{T}) \stackrel{f^*} \rightarrow
Sym(\hat{S})$,
or in other words, a map $f^*: R_S \rightarrow R_T$.

\begin{lemma} \label{t-map}(cf. \cite{A-B})
Suppose there is a $T$-map
$X \stackrel{\phi} \rightarrow S$.
Then $t\cdot A_*^T(X,m)= 0$ for any $t=f^*s$ with $s \in R_S^+$.
\end{lemma}

Proof of Lemma \ref{t-map}.
Since $A^*_S$ is generated in degree 1, we may
assume that $s$ has degree 1. After clearing denominators
we may assume that $s = c_1(L_s)$ for some line bundle
on a space $U/S$. The action of  $t=f^*s$ on $A_*(X_T)$ is just
given by $c_1(\pi_T^*f^*L_s)$ where $\pi_T$ is the
map $U \times^T X \rightarrow U/T$.
To prove the lemma we will show that this bundle is trivial.

First note that  $L_s = U \times^S k$ for some action of $S$ on
the one-dimensional vector space $k$.
The pullback bundle on $X_T$ is the line bundle
$$U \times^T(X \times k) \rightarrow X_T$$
where $T$ acts on $k$ by the composition of $f:T \rightarrow S$
with the original $S$-action.
Now define a map
$$\Phi: X_T \times k \rightarrow U \times^T(X \times k)$$
by the formula
$$\Phi(e,x,v) = (e,x,\phi(x)\cdot v)$$
(where $\phi(x) \cdot v$ indicates the original $S$ action).
This map is well defined since
\begin{eqnarray*}
\Phi(et,t^{-1}x,v) & = &(et, t^{-1}x, \phi(t^{-1}x) \cdot v)\\
& = & (et,t^{-1}x,t^{-1} \cdot(\phi(x) \cdot v))
\end{eqnarray*}
as required. This map is easily seen to be an isomorphism
with inverse $(e,x,v) \mapsto (e,x,\phi(x)^{-1} \cdot v)$.
\endproof

\begin{prop} \label{fix}
If $T$ acts on $X$
without fixed points, then there exists $r \in R_T^+$ such that
$r \cdot A_*^T(X,m)= 0$. (Recall that $A_*^T(X,m)$ refers to
$T$-equivariant higher Chow groups.)
\end{prop}

\medskip

Before we prove Proposition \ref{fix}, we state and prove
a lemma.

\begin{lemma} \label{porb}
If $X$ is a variety with an action of a torus $T$, then there is
an open $U \subset X$ so that the stabilizer is constant for
all points of $U$.
\end{lemma}

Proof of Lemma \ref{porb}. Let $\tilde X \rightarrow X$ be the normalization
map. This map is $T$-equivariant and is an isomorphism
over an open set. Thus we may assume $X$ is normal. By Sumihiro's
theorem, the $T$ action on $X$ is locally linearizable, so it
suffices to prove the lemma when $X = V$ is a vector space and
the action is diagonal.

If $V = k^n$, then let $U = (k^*)^n$. The $n$-dimensional
torus ${\bf G}_m^n$ acts transitively on $U$ in the obvious way.
This action commutes with the given action of $T$. Thus the stabilizer
at each closed point of $U$ is the same.
\endproof

Proof of Proposition \ref{fix}. Since $A_*^G(X) = A_*^G(X_{red})$
we may assume $X$ is reduced. Working with each component
separately, we may assume $X$ is a variety. Let $X^0 \subset X$
be the ($G$-invariant) locus of smooth points.
By Sumhiro's theorem \cite{Sumihiro}, the action of
a torus on a normal variety is locally
linearizable (i.e. every point has an affine invariant
neighborhood). Using this theorem it is easy to see
that the set
$X(T_1) \subset X^0$ of points with stabilizer
$T_1$ can be given the structure of a locally closed
subscheme of $X$.
By Lemma \ref{porb} there is
some $T_1$ such that $U= X(T_1)$ is open in $X^0$, and thus in $X$.

The torus $T'=T/T_1$ acts without stabilizers,
but the action of $T'$ on $U$ is not a priori proper. However,
by \cite[Proposition 4.10]{Th1}, we can replace $U$
by a sufficiently small open set so that $T'$ acts freely
on $U$ and a principal bundle quotient $U \rightarrow U/T$
exists. Shrinking $U$ further, we can assume that this
bundle is trivial, so there is a $T$ map $U \rightarrow T'$.
Hence, by the lemma, $t \cdot A^T_*(U) = 0$ for any $t \in A_T^*$
which is pulled back
from $A^*_{T'}$.

Let $Z = X -U$.  By induction on dimension, we may assume $p \cdot
A^T_*Z = 0$ for some homogeneous polynomial $p \in R_T$.  From the
long exact sequence of higher Chow groups,
$$\ldots A_*^T(Z,m) \rightarrow A_*^T(X,m) \rightarrow A_*^T(U,m) \rightarrow
\ldots$$
it follows that $tp$ annihilates $A_*^T(X)$ where
$t$ is the pullback of a homogeneous element of degree $1$
in $R_S$.
\endproof

If $X$ is a scheme with a $T$-action,
we may put a closed subscheme structure
on the locus $X^T$ of points fixed by
$T$ (\cite{Iversen}).
Now $R_T= Sym(\hat{T})$ is a polynomial ring.
Set ${\cal Q}= (R_T^+)^{-1} \cdot R_T$, where $R_T^+$
is the multiplicative system of homogeneous elements of positive degree.
\begin{thm} \label{lcztn}(localization)
The map $i_*:A_*(X^T)  \otimes {\cal Q} \rightarrow A_*^T(X) \otimes
{\cal Q}$ is
an isomorphism.
\end{thm}

\medskip

Proof of Theorem \ref{lcztn}.
By Proposition \ref{fix}, $A_*^T(X - X^T,m) \otimes {\cal Q} = 0$.
Thus by the localization exact sequence
$A_*^T(X^T) \otimes {\cal Q} = A_*^T(X) \otimes {\cal Q}$
as desired. \endproof

{\bf Remark.}
The strategy of the
proof is similar to proofs in other theories, see for example
for \cite[Chapter 3.2]{Hsiang}.

\section{Explicit localization for smooth varieties}
The localization theorem in equivariant cohomology has a more explicit
version for smooth varieties because the
fixed locus is regularly embedded.
This yields an integration formula from which
the Bott residue formula is easily deduced (\cite{A-B}, \cite{B-V}).
In this section we prove the analogous results for equivariant Chow
groups of smooth varieties.  Because equivariant Chow theory has
formal properties similar to equivariant cohomology, the arguments are
almost the same as in \cite{A-B}. As before we assume that all tori
are split.

Let $F$ be a scheme with a trivial $T$-action.
If $E \rightarrow F$ is a $T$-equivariant vector bundle on
$F$, then $E$ splits canonically into a direct sum of vector subbundles
$\oplus_{\lambda \in \hat{T}} E_{\lambda}$, where $E_{\lambda}$
consists of the subbundle of vectors in $E$ on which $T$ acts by the
character $\lambda$.  The equivariant Chern classes of an
eigenbundle $E_{\lambda}$ are given by the following lemma.

\begin{lemma} \label{l.trivchern}
Let $F$ be a scheme with a trivial $T$-action, and let
$E_{\lambda} \rightarrow F$ be a $T$-equivariant vector bundle of rank
$r$ such
that the action of $T$ on each vector in $E_{\lambda}$ is given by the
character $\lambda$.  Then for any $i$,
$$
c^T_i(E_{\lambda}) = \sum_{j \leq i}
\left( \begin{array}{c} r-j \\
i-j \end{array}
\right)
c_j(E_{\lambda}) \lambda^{i-j}.
$$
In particular the component of $c_r^T(E_{\lambda})$ in $R^r_T$ is
given by $\lambda^r$.  \endproof
\end{lemma}

As noted above, $A^*_T(F) \supset A^*F \otimes R_T$.  The lemma
implies that $c^T_i(E)$ lies in the subring $A^*F \otimes R_T$.
Because $A^N F = 0$ for $N > \mbox{dim }F$, elements of $A^i F$, for
$i>0$, are nilpotent elements in the ring $A^*_T(F)$.  Hence an
element $\alpha \in A^d F \otimes R_T$ is invertible in $A^*_T(F)$ if
its component in $A^0 F \otimes R^d_T \cong R^d_T$ is nonzero.

For the remainder of this section $X$ will denote a smooth variety
with a $T$ action.
\begin{lemma} \cite{Iversen}
If $X$ is smooth then the fixed locus
$X^T$ is also smooth. \endproof
\end{lemma}
For each component $F$ of the fixed locus $X^T$
the normal bundle $N_FX$ is a $T$-equivariant vector bundle over $F$.
Note that the action of $T$ on $N_FX$ is non-trivial.

\begin{prop}
If $F$ is a component of $X^T$
with codimension $d$ then $c_d^T(N_FX)$ is invertible
in $A^*_T(F) \otimes {\cal Q}$.
\end{prop}
Proof: By (\cite[Proof of Proposition 1.3]{Iversen}),
for each closed point $f \in F$, the tangent space
$T_fF$ is equal to $(T_fX)^T$, so $T$ acts with non-zero weights on the
normal space $N_f = T_fX/T_fF$. Hence the characters $\lambda_i$
occurring in the
eigenbundle decomposition of $N_FX$ are all non-zero.  By the
preceding lemma, the component of $c_d^T(N_FX)$ in $R^d_T$ is nonzero.
Hence $c_d^T(N_FX)$ is invertible
in $A^*_T(F) \otimes {\cal Q}$, as desired.  \endproof

\medskip

Using this result we can get, for $X$ smooth,
the following more explicit version
of the localization theorem.

\begin{thm} \label{xxx}(Explicit localization)
Let $X$ be a smooth variety with a torus
action.
Let $\alpha \in A_*^T(X) \otimes {\cal Q}$.
Then $$\alpha = \sum_F
i_{F*}\frac{i^*_F\alpha}{c_{d_F}^T(N_FX)},$$ where the sum is over the
components $F$ of $X^T$ and $d_F$ is the codimension of $F$ in $X$.
\end{thm}
Proof: By the surjectivity part of the localization theorem,
we can write $\alpha = \sum_F
i_{F*}(\beta_F)$.  Therefore, $i^*_F\alpha = i^*_Fi_{F*}(\beta_F)$
(the other components of $X^T$ do not contribute); by the
self-intersection formula, this is equal to $ c_{d_F}^T(N_FX) \cdot
\beta_F$.  Hence $\beta_F = \frac{i^*_F\alpha}{c_{d_F}^T(N_FX)}$ as
desired. \endproof

{\bf Remark.} This formula is valid, using the virtual normal bundle,
even if $X$ is singular, provided that the embedding of the fixed
locus in $X$ is a local complete intersection morphism. Unfortunately,
this condition is difficult to verify.  However, if $X$ is cut by a
regular sequence in a smooth variety, and the fixed points are
isolated, then the methods of \cite[Section 3]{BFQ} can be used to
give an explicit localization formula.  A similar remark applies to
the Bott residue formula below.

\medskip

If $X$ is complete, then the projection $\pi_X: X \rightarrow pt$
induces push-forward maps $\pi_{X*}: A^T_* X \rightarrow R_T$ and
$\pi_{X*}: A^T_* X \otimes {\cal Q} \rightarrow {\cal Q}$.  There
are similar maps with $X$ replaced with any component $F$ of $X^T$.
Applying $\pi_{X*}$ to both sides of the explicit localization
theorem, and noting that $\pi_{X*} i_{F*} = \pi_{F*} $, we deduce
the ``integration formula'' (cf. \cite[Equation (3.8)]{A-B}).

\begin{cor}
(Integration formula) Let $X$ be smooth and complete, and
let $\alpha \in A_*^T(X) \otimes {\cal Q}$.  Then
$$\pi_{X*}(\alpha) = \sum_{F \subset X^T} \pi_{F*}\left(
\frac{i^*_F\alpha}{c_{d_F}^T
(N_FX)}\right)$$
as elements of ${\cal Q}$.  \endproof
\end{cor}

\medskip

{\bf Remark.} If $\alpha$ is in the image of the natural map $A_*^T(X)
\rightarrow A_*^T(X) \otimes {\cal Q}$ (which need not be injective),
then the equation above holds in the subring $R_T$ of ${\cal Q}$.  The
reason is that the left side actually
lies in the subring $R_T$; hence so does the right side.  In the
results that follow, we will have expressions of the form $z = \sum
z_j$, where the $z_j$ are degree zero elements of ${\cal Q}$ whose sum
$z$ lies in the subring $R_T$.  The pullback map from equivariant to
ordinary Chow groups gives a map $i^*: R_T = A^T_* (pt) \rightarrow \Q
= A_* (pt)$, which identifies the degree 0 part of $R_T$ with $\Q$.
Since $\sum z_j$ is a degree 0 element of $R_T$, it is identified via
$i^*$ with a rational number.  Note that $i^*$ cannot be applied to
each $z_j$ separately, but only to their sum.  In the integration and
residue formulas below we will identify the degree 0 part of $R_T$
with $\Q$ and suppress the map $i^*$.  \medskip

The preceding corollary yields an integration formula for
an element $a$ of the ordinary Chow group $A_0 X$, provided that $a$ is
the pullback of an element $\alpha \in A^T_0 X$.

\begin{prop}
Let $a \in A_0 X$, and suppose that $a = i^* \alpha$ for $\alpha \in
A^T_0 X$.  Then
$$
\mbox{deg }(a) = \sum_F  \pi_{F*}\{\frac{i^*_F\alpha}{c_{d_F}^T
(N_FX)} \}
$$
\end{prop}
Proof: Consider the commutative diagram
$$\begin{array}{ccc}
X & \stackrel{i} \hookrightarrow & X_T\\
\downarrow\scriptsize{\pi_X} & & \downarrow\scriptsize{\pi^T_X}\\
\mbox{pt} & \stackrel{i} \rightarrow & U/T .
\end{array}$$
We have $\pi_{X*}(a) = \pi_{X*} i^*(\alpha) = i^* \pi^T_{X*}(\alpha)$.
Applying the integration formula gives the result. \endproof

\subsection{The Bott residue formula}
Let $E_1, \ldots , E_s$ be a $T$-equivariant vector bundles
on a complete, smooth $n$-dimensional variety $X$.
Let $p(x^1_1, \ldots x^1_s,\ldots , x^n_1, \ldots x^n_s)$
be a polynomial of weighted degree $n$,
where $x^i_j$ has weighted degree $i$.
Let $p(E_1, \ldots , E_s)$ denote the polynomial
in the Chern classes of $E_1, \ldots , E_s$ obtained
setting $x^i_j = c_i(E_g)$.
The integration formula above will allow us to compute
$\mbox{deg }(p(E_1, \ldots , E_s) \cap [X])$
in terms of the
restriction of the $E_i$ to $X^T$.

As a notational shorthand, write $p(E)$
for $p(E_1, \ldots, E_s)$ and $p^T(E)$ for the corresponding
polynomial in the $T$-equivariant Chern classes of $E_1, \ldots
, E_r$.
Notice that
$p(E) \cap [X] = i^* (p^T(E) \cap [X_T])$.
We can therefore apply the
preceding proposition to get the Bott residue formula.

\begin{thm} \label{bott}
(Bott residue formula) Let $E_1, \ldots , E_r$ be a $T$-equivariant
vector bundles a complete, smooth $n$-dimensional
variety. Then
$$
\mbox{deg }(p(E) \cap [X]) =  \sum_{F \subset X^T}
\pi_{F*}\left(\frac{p^T(E|_{F}) \cap
[F]_T}{c_{d_F}^T (N_FX)} \right).
$$
\endproof
\end{thm}

{\bf Remark.} Using techniques of
algebraic deRham homology, H\"ubl and Yekutieli \cite{H-Y} proved a
version of the Bott residue formula, in characteristic 0,
for the action of any algebraic vector
field with isolated fixed points.

\medskip

By Lemma \ref{l.trivchern} the
equivariant Chern classes $c^T_i(E_j|_{F})$ and $c_{d_F}^T (N_FX)$ can
be computed in terms of the characters of the
torus occurring in the eigenbundle
decompositions of $E_j|_{F}$ and $N_FX$ and the Chern classes of the
eigenbundles.
The above formula can then be readily converted (cf. \cite{A-B}) to more
familiar forms of the Bott residue formula not involving equivariant
cohomology. We omit the details. If the torus $T$ is
1-dimensional, then degree zero elements of ${\cal Q}$ are rational
numbers, and the right hand side of the formula is just a sum of
rational numbers. This is the form of the Bott residue formula which
is most familiar in practice.

\section{Localization and residue formulas for singular varieties}
\label{singex}

In general, the problem of proving localization and residue formulas
on singular varieties seems interesting and difficult.  In this
section we discuss what can be deduced from an equivariant embedding
of a singular scheme $X$ into a smooth $M$.  The results are not very
general, but (as we show) they can be applied in some interesting
examples, for example, if $X$ is a complete intersection in $M = \P^n$
and $T$ acts on $M$ with isolated fixed points, or if $X$ is a
Schubert variety in $M = G/B$.

The idea of using an embedding into a smooth variety to extract
localization information is an old one. In the case of the action of
an automorphism of finite order, the localization and Lefschetz
Riemann-Roch formulas of \cite{Quart}, \cite{BFQ} on quasi-projective
varieties are obtained by a calculation on $\P^n$. Moreover, as in our
case, the best formulas on singular varieties are obtained when the
embedding into a smooth variety is well understood.

At least in principle, a localization theorem can be deduced if every
component of $X^T$ is a component of $M^T$.  This holds, for example,
if the action of $T$ on $M$ has isolated fixed points; or if $X$ is a
toric (resp.  spherical) subvariety of a nonsingular toric
(resp. spherical) variety $M$.  In particular, the condition holds if
$X$ is a Schubert variety and $M$ is the flag variety.  We have the
following proposition.

\begin{prop} \label{xxxsing}
Let $f: X \rightarrow M$ be an equivariant embedding of $X$ in a nonsingular
variety $M$.  Assume that every component of $M^T$ which intersects $X$ is
contained in $X$.  If $F$ is a component of $X^T$, write $i_F$ for the
embedding of $F$ in $X$, and $j_F$ for the embedding of $F$ in $M$.  Then:

$(1)$ $f_*: A_*^T(X) \otimes {\cal Q} \rightarrow A_*^T(M) \otimes {\cal Q}$
is injective.

$(2)$ Let $\alpha \in A_*^T(X) \otimes {\cal Q}$.
Then $$\alpha = \sum_F
i_{F*}\frac{j^*_F f_* \alpha}{c_{d_F}^T(N_FM)},$$ where the sum is over the
components $F$ of $X^T$ and $d_F$ is the codimension of $F$ in $M$.
\end{prop}

Proof: (1) Since the components of $X^T$ are a subset of
the components of $M^T$,
$\oplus_{F \subset X^T} A_*^T(F)$ is an $R_T$-submodule of $\oplus
_{F \subset M^T} A_*^T(F)$. By the localization theorem,
$$\sum_{F \subset X^T} j_{F*}(A_*^T(F)) \otimes {\cal Q} \simeq A_*^T(X)
\otimes {{\cal Q}}$$
and
$$\sum_{F \subset M^T} i_{F*}(A_*^T(F)) \otimes {\cal Q} \simeq A_*^T(M)
\otimes {{\cal Q}}.$$
Since $i_{F*} = f_* i_{F *}$, the result follows.

(2) By (1) it suffices to prove
that $$f_*(\alpha - \sum_F
j_{F*}\frac{i^*_F f_* \alpha}{c_{d_F}^T(N_FM)}) = 0 \in A_*^T(X) \otimes
{\cal Q}.$$
Since $i_{F*} = f_* j_{F*}$ the theorem follows from the explicit
localization theorem applied to the class $f_* \alpha$
on the smooth variety $M$.
\endproof

To obtain a residue formula that computes
Chern numbers of bundles on $X$, we only need to know an expansion
$[X]_T =\sum_{F\subset X^T} i_{F*}(\beta_F)$,
where $\beta_F \in A_*^T(F)$. In this case we obtain the formula
$$
\mbox{deg }(p(E) \cap [X]) =  \sum_{F \subset X^T}
\pi_{F*}\left(\frac{p^T(i_F^*(E)) \cap \beta_F}{c_{d_F}^T(N_FM)}\right).$$

In the setting of Proposition \ref{xxxsing}, the classes $\beta_F$ are
given by $i_F^* f_* [X]_T$.  To obtain a useful residue formula, we
need to make this expression more explicit.  This is most easily done
if we can express $f_* [X]_T$ in terms of Chern classes of naturally
occuring equivariant bundles on $M$.  The reason is that the pullback
$i_F^*$ of such Chern classes is often easy to compute, particularly
if $F$ is an isolated fixed point (cf. Lemma. \ref {l.trivchern}).
Indeed, this is why the Bott residue formula is a good calculational
tool in the non-singular case.

%
%

Although the conditions to obtain localization and residue formulas
are rather strong, they are satisfied in some interesting cases.  We
will consider in detail two examples: complete intersections in
projective spaces, and Schubert varieties in $G/B$. For complete
intersections some intrinsic formulas can be deduced using the virtual
normal bundle (see the remark after Theorem \ref{xxx}).  In this
section our point of view for complete intersections is different. We
do not use the virtual normal bundle, but instead use the fact that if
$X \stackrel{f} \hookrightarrow M$ is a complete intersection, it is
easy to calculate $f_*[X]_T \in A_*^T(M)$.  As an example of our methods
we do a localization and residue calculation on a singular quadric in
$\P^3$.

As a final remark, note that to compute Chern numbers of bundles on
$X$ which are pulled back from $M$, it suffices to know $f_*[X]_T$,
for then we can apply residue formulas on $M$.  Information about the
fixed locus in $X$ is irrelevant.  The interesting case is when the
bundles are not pulled back from $M$; see the example of the singular
quadric below.

\subsection{Complete intersections in projective space}

For simplicity we consider the case where the dimension of $T$ is $1$.
If $T$ acts on a vector space $V$ with weights $a_0, \ldots, a_n$ then
$A^*_T(\P(V)) = \Z[h,t] / \prod (h + a_i t)$.  We are interested in
complete intersections $X$ in $\P(V)$ where the functions $f_i$
defining $X$ are, up to scalars, preserved by the $T$-action, i.e.,
$t \cdot f_i = t^{a_i} f_i$.  In this case we say $f_i$ has weight $a_i$.

The following lemma is immediate.

\begin{lemma}
Suppose $X$ is a hypersurface in $\P(V)$ defined by a homogeneous
polynomial $f$ of degree $d$ and weight
$a$.  Then $[X]_T = d h + a t \in A^*_T(\P(V))$.  Hence if $X$ is a
complete intersection in $\P(V)$ defined by homogeneous polynomials
$f_i$ of degree $d_i$ and weight $a_i$,
then $[X]_T = \prod (d_i h + a_i t)$. \endproof
\end{lemma}

If $T$ acts on $V$ with distinct weights, then $T$ has isolated fixed
points on $M = \P(V)$, and (trivially) every component of $X_T$ is a
component of $M^T$; so by the preceding discussion there is a useful
residue formula. In particular using a little linear algebra we can
easily obtain a formula for $[X]_T$ in terms of the fixed points in
$X$ and the weights of the action. We omit the details to avoid a
notational quagmire, but the ideas are illustrated in the example
of the singular quadric.

\subsection{Schubert varieties in $G/B$} \label{schubs}

In this section, we work over an algebraically closed field.  For
simplicity, we take Chow groups to have rational coefficients, and
let $R = R_T \otimes \Q$ denote the rational equivariant Chow ring
of a point.

Let $G$ be a reductive group and $B$ a Borel subgroup, and $\B = G/B$
the flag variety.  In the discussion below, the smooth variety $\B$ will
play the role of $M$, and the Schubert variety $X_w$ the role of $X$.

Let $T \subset B$ be a maximal torus.  $T$ acts on $\B$ with finitely
many fixed points, indexed by $w \in W$; denote the corresponding
point by $p_w$.  More precisely, if we let $w$ denote both an element
of the Weyl group $W = N(T)/T$ and a representative in $N(T)$, then
$p_w$ is the coset $wB$.  The flag variety is a disjoint union of the
$B$-orbits $X_w^0 = B \cdot p_w$.  The $B$-orbit $X_w^0$ is called a
Schubert cell and its closure $X_w$ a Schubert variety.  If $e$
denotes the identity in $W$ and $w_0$ the longest element of $W$, then
$X_e$ is a point and $X_{w_0} = \B$.

We have $X_w = \cup_{u \leq w} X_u^0$.  The $T$ equivariant Chow group
of $X_w$ is a free $R_T$-module with basis $[X_u]_T$, for $u \leq w$.

Let $j_u: p_u \hookrightarrow \B$.  Fix $w \in W$ and let $f: X_w
\hookrightarrow \B$.  For $u \leq w$ let $i_u: p_u \hookrightarrow X_u$.
If $v \leq w$ let $[X_v]_T$ denote the equivariant fundamental class
of $X_v$ in $A_*^T(X_w)$.  We want to make explicit the localization
theorem for the variety $X_w$ (which is singular in general), i.e., to
compute $[X_v]_T$ in terms of classes $i_{u*} \beta_u$.

The (rational) equivariant Chow groups $A^*_T(\B)$ can be described as
follows.  We consider two maps $\rho_1, \rho_2: R \rightarrow
A^*_T(\B)$.  The map $\rho_1$ is the usual map $R \rightarrow
A^*_T(\B)$ given by equivariant pullback from a point.  The definition
of $\rho_2$ is as follows.  For each character $\lambda \in \hat{T}$
set $\rho_2(\lambda) = c_1^T(M_\lambda)$ where $M_{\lambda}$ is the
line bundle $G \times^B k_{\lambda} \rightarrow \B$; extend $\rho_2$
to an algebra map $R \rightarrow A^*_T(\B)$.  The map $R \otimes_{R^W}
R \rightarrow A^*_T(\B)$ taking $r_1 \otimes r_2$ to $\rho_1(r_1)
\rho_2(r_2)$ is an isomorphism (see e.g. \cite{Brion}).

We adopt the convention that the Lie algebra of $B$ contains the
positive root vectors.  We can identify $T_{p_w}(\B)$ with $\g /
(\Ad w) \b$.  This is a representation of $T$ corresponding to the
$T$-equivariant normal bundle of the fixed point $p_w$.  We identify
$A^*_T(p_w) \cong R$.  From our description of $T_{p_w}(\B)$, we
see that (if $n$ denotes the dimension of $\B$) $c_n^T(N_{p_w}\B)$ is
the product of the roots in $\g /(\Ad w) \b$, which is easily
seen to give
$$
c_n^T(N_{p_w}\B)  = c_w := (-1)^n(-1)^w \prod_{\alpha > 0} \alpha.
$$
where $n$ is the number of roots $\alpha >0$.

To obtain a localization formula we also need to know the maps $j_u^* :
A^*_T(\B) \rightarrow A^*_T(p_u)$, where $j_u: p_u \hookrightarrow \B$
is the inclusion.  We have identified $A^*_T(\B) = R \otimes_S R$ and
$A^*_T(p_w) = R$.  Thus, we may view $j_u^*$ as a map $R \otimes_S R
\rightarrow R$.  There is a natural action of $W \times W$ on $R
\otimes_S R$.  Let $m: R \otimes_S R \rightarrow R$ denote the
multiplication map.

\begin{lemma}
For $u \in W$, the map $j_u^*: R \otimes_S R \rightarrow R$ equals the
composition $m \circ (1 \times u)$.
\end{lemma}

Proof: It suffices to show that $j_u^* \rho_1(\lambda) = \lambda$ and
$j_u^* \rho_2(\lambda) = u \lambda$.  Now, $j_u^* \rho_1$ is just the
equivariant pullback by the map $p_u \rightarrow pt$.  Since this
equivariant pullback is how we identify $A^*_T(p_u) = A^*_T(pt) = R$,
with these identifications, $j_u^*\rho_1$ is the identity map, $j_u^*
\rho_1(\lambda) = \lambda$.  Also, by definition $j_u^*
\rho_2(\lambda) = c_1^T(M_{\lambda}|_{p_u})$.  As a representation of
$T$, $M_{\lambda}|_{p_u} \cong k_{u \lambda}$, so $c_1^T(M_{\lambda}|_{p_u})
= u \lambda$, as desired.  \endproof

If $F \in R \otimes_S R$ is a polynomial set
$F(u) = j_u^*F \in R$.

Recall that we have fixed $w$ and let $f : X_w \rightarrow \B$ denote
the inclusion; for $v \leq w$, $[X_v]_T$ denotes a class in
$A_*^T(X_w)$.  By work of Fulton and Pragacz-Ratajski, for $G$
classical, it is known how to express $f_*[X_v]_T \in A^*_T(\B)$ in
terms of the isomorphism $A^*_T(\B) \cong R \otimes_{R^W} R$.  More
precisely, Fulton and Pragacz-Ratajski (\cite{Fu1}, \cite{P-R})
define elements in $R
\otimes_{\Q} R$ which project to $f_*[X_v]_T$ in $R \otimes_{R^W} R$.
Let $F_u$ denote either the polynomial defined by Fulton or that
defined by Pragacz-Ratajski.  Using these polynomials we can get an
explicit localization formula for Schubert varieties.

\begin{prop}
With notation as above, the class $[X_v]_T$ in $A_*^T(X_w) \otimes
{\cal Q}$ is given by
$$
[X_v]_T = (-1)^n
\frac{1}{\prod_{\alpha > 0} \alpha} \sum_{u \leq v} (-1)^u i_{u*}
\left( F_v(u) \cap [p_u]_T \right).
$$
\end{prop}

Proof: This is an immediate consequence of the preceding discussion and
Proposition \ref{xxxsing}.  \endproof

{\bf Remark.} Taking $w = w_0$, so $X_w = \B$, the above formula is an
explicit inverse to the formula of \cite[Section 6.5, Proposition
(ii)]{Brion}.  This shows $\frac{f_w(u)}{\Pi_{\alpha > 0} \alpha}$
is Brion's equivariant multiplicity of $X_w$ at the fixed point $p_u$, and
also links Brion's proposition to \cite[Theorem 1.1]{G}.

\subsection{A singular quadric}

In this section we consider the example of the singular quadric $Q
\stackrel{f} \hookrightarrow \P^3$ defined by
the equation $x_0x_1 + x_2^2 = 0$ (note that we allow the
characteristic to be 2).  Let $\P^2 \subset \P^3$ be the hyperplane
defined the equation $x_2 = 0$ and let $\pi: Q \rightarrow \P^2$ be
the projection from $(0,0,1,0)$.  As a sample of the kinds of the
residue calculations that are possible, we prove the following
proposition.

\begin{prop}
 $$\int_Q c_1(\pi^* T_{\P^2}) c_1(f^*T_{\P^3}) = 24$$.
\end{prop}
Proof. We will prove this by considering the following torus action.
Let $T = \G_m$ act on
$\P^3$ with weights $(1,-1,0,a)$, where $a \notin \{0,-1,1\}$.  The
quadric is invariant under this action, so $T$ acts on $Q$.
Since $(0,0,1,0)$ is a fixed point, $\pi$ is an equivariant map
where $T$ acts on $\P^2$ with weights $(1,-1,a)$.
Thus $c_1^T(\pi^*T_{\P^2}) c_1^T(f^*T_{\P^3}) \cap [Q]_T$
defines an element of $A_*^T(Q) \otimes {\cal Q}$ which
we will express as a residue in terms of the fixed points
for the action of $T$ on $Q$. To do this
we need to express $[Q]_T$ in terms of the fixed points.
By Proposition \ref{xxxsing} this can be done if we know
$f_*[Q]_T \in A_*^T(\P^3)$. Since $Q$ is a quadric
of weighted degree 0 with respect to the $T$-action,
$f_*[Q]_T = 2h \in A_*^T(\P^3)$.

Since everything can be done explicitly, we will calculate
more than we need and determine the entire $R_T$-module $A_*^T(Q)$
in terms of the fixed points.

\medskip

{\bf Explicit localization on the singular quadric.}
The quadric has a decomposition into affine cells with one cell in
dimensions 0,1 and 2.  The open cell is $Q_0 = \{(1,-x^2,x,y) | (x,y)
\in k^2\}$.  In dimension 1 the cell is $l_0 = \{(0,1,0,x) | x \in
k\}$, and in dimension 0, the cell is the singular point
\footnote{In characteristic 2, this point is not an isolated singular
point.} $p_s =
\{(0,0,0,1)\}$.  Thus, $A_i(Q) = \Z$ for $i = 0,1,2$ with generators
$[Q]$, $[l]$ and $[p_s]$. Moreover, these cells are $T$-invariant, so
their equivariant fundamental classes form a basis for $A_*^T(Q)$ as
an $R_T = \Z[t]$ module.  Let $\1$, $L$, and $P_s$ denote the
corresponding equivariant fundamental classes $[Q]_T,[l]_T$, and
$[p_s]_T$.

There are three fixed points $p_s = (0,0,0,1)$, $p = (1,0,0,0)$ and
$p' = (0,1,0,0)$.  These points have equivariant fundamental classes
in $A_*^T(Q^T)$ which we denote by $P_s$, $P$, and $P'$.
By abuse of notation we will not distinguish between
$P_s$ and $i_*(P_s)$.

Both $A_*^T(Q^T)$ and $A_*^T(Q)$ are free $R_T$-modules of rank $3$,
with respective ordered bases $\{P_s,P,P' \}$ and $\{
P_s,L,\1 \}$.  The map $i_*$ is a linear transformation of these
$R_T$-modules, and we will compute its matrix with respect to these
ordered bases.  This matrix can be easily inverted, provided we invert
$t$, and so we obtain $(i_*)^{-1}$.

The equivariant Chow ring of $\P^3$ is given by
$$
A^*_T(\P^3) = \Z[t,h]/(h-t)(h+t)h(h+at)
$$
so $A^*_T(\P^3)$ is free of rank $4$ over $R_T$, with basis
$\{1,h,h^2,h^3\}$.

To compute $i_*P_s$, $i_*P$, and $i_*P'$, we take advantage of the
fact that the pushforward $f_*: A_*^T(Q) \rightarrow A_*^T(\P^3)$ is
injective.  Moreover it is straightforward to calculate the
pushforward to $\P^3$ of all the classes in our story.  To simplify
the notation, we will use $f_*$ to denote either of the maps $A_*^T(Q)
\rightarrow A_*^T(\P^3)$ or $A_*^T(Q^T) \rightarrow A_*^T(\P^3)$.  We find:
$$\begin{array}{l}
f_*(\1) = 2h\\
f_*(L)= (h-t)h\\
f_*(P_s) = h^3 - ht^2\\
f_*(P) = h^3 + (a-1) h^2 t - a h t^2\\
f_*(P') = h^3 + (a+1) h^2 t + a h t^2
\end{array} $$
This implies that
$$\begin{array}{l}
i_*(P_s) = P_s\\
i_*(P) = (a-1)t L + P_s\\
i_*(P') = (a+1)t^2 \1 + (a+1)t L + P_s.
\end{array} $$
So the matrix for $i_*^T$ is
$$\left(\begin{array}{ccc}  1 & 1 & 1\\ 0 & (a-1)t  & (a+1)t\\ 0 & 0 &
(a + 1)t^2 \end{array} \right).$$
Inverting this matrix we obtain
$$\left(\begin{array}{ccc} 1 & \frac{1}{t(- a)} & \frac{2}{t^2(a^2 -1)} \\
0 & \frac{1}{t(a-1)} & \frac{1}{t^2(1-a)} \\
0 & 0 & \frac{1}{t^2(a+1)} \end{array} \right)
$$
Thus we can write (after supressing the $(i_*)^{-1}$ notation)
$$\begin{array}{l}
P_s = P_s\\
L = \frac{1}{t(a-1)}( -P_s + P)\\
\1 =\frac{1}{t^2(a^2-1)}(2P_s - (a+1)P  +(a -1)P')
\end{array}. $$

\medskip
{\bf Calculation of Chern numbers}

We now return to the task of computing $c_1(\pi^*T_{\P^2}) c_1(f^*T_{\P^3})
\cap \1$. To simplify notation, set $\alpha_1 := c_1(\pi^*T_{\P^2})$
and $\alpha_2 := c_1(f^*T_{\P^3})$ and $\alpha := \alpha_1 \alpha_2$.
By the calculations above
$$ \alpha_1 \alpha_2 \cap \1 = i_*(i^*\alpha_1 i^*\alpha_1 \cap
\frac{t^{-2}}{a^2-1}(2P_s - (a+1)P  +(a -1)P')).$$
To compute the class explicitly we must compute the restrictions
of $\alpha_1$ and $\alpha_2$ to each of the fixed points $P_s$, $P$ and
$P'$.

The tangent space to $P_s$ in $\P^3$ has weights $(1-a, -1-a, -a)$.
Thus $\alpha_{2}|_{P_s} = (1-a)t  - (1+a)t - at = -3at$.
To compute $\alpha_{1}|_{P_s}$ observe that $P_s$ is the
inverse image of the fixed point $(0,0,1) \in \P^2$. Since
$T_{\P^2}$ has weights $(1-a,-1-a)$ at this point,
$c_1(\pi^*T_{\P_2})|_{P_s} = (1-a)t + (-1 -a)t = -2t$.

The restrictions to the other two fixed points can be calculated
similarly. In particular
$$\begin{array}{ll}
\alpha_{1}|_{P} = (a-3)t \mbox{  }& \alpha_{1}|_{P'} = (a+3)t\\
\alpha_{2}|_{P} = (a- 4)t  & \alpha_{2}|_{P'} = (a+4)t
\end{array}
$$
Thus,
$$\alpha \cap \1 = \frac{12 a^2}{a^2-1}P_s - \frac{(a-3)(a-4)(a+1)}{a^2-1}P
+ \frac{(a+3)(a+4)(a-1)}{a^2-1}P' \in A_*^T(Q) \otimes \Q.$$
Thus,
$$\begin{array}{ll}
\mbox{deg }(c_1(\pi^*T_{\P^2}) c_1(f^{*}T_{\P^3}) \cap [Q]) & =
\frac{12 a^2}{ a^2 -1 } \;- \;\frac{(a-3)(a-4)(a+1)}{a^2-1} \; + \;
\frac{(a+3)(a+4)(a-1)}{a^2-1}\\
& = 24
\end{array}. $$ \endproof

\end{document}